\documentclass[aos,MSNbibl,citesort,dvips]{arximspdf}

\doi{10.1214/11-AOS945}
\volume{39}
\issue{6}
\pubyear{2011}
\firstpage{3357}
\lastpage{3368}

\makeatletter

\newtheorem{lemma}{Lemma}

\newtheorem{theorem}{Theorem}

\newproclaim{Remark}{Remark}

\newcommand{\half}{\frac{1}{2}}
\newcommand{\bmz}{{\mathbf0}}
\newcommand{\bx}{{\mathbf x}}
\newcommand{\bu}{{\mathbf u}}
\newcommand{\bX}{{\mathbf X}}
\newcommand{\be}{{\mathbf e}}
\newcommand{\byr}{{\mathbf y}}
\newcommand{\bz}{{\mathbf z}}
\newcommand{\bI}{{\mathbf I}}
\newcommand{\bj}{{\mathbf j}}
\newcommand{\ba}{{\mathbf a}}
\newcommand{\bS}{{\mathbf S}}
\newcommand{\bs}{{\mathbf s}}
\newcommand{\bmg}{{\bolds\mu}}
\newcommand{\Prob}{\mathrm{P}}

\newcommand{\cgfargv}{{\bolds\tau}}
\newcommand{\Vol}{\operatorname{Vol}}
\newcommand{\reals}{\mathbb R}
\newcommand{\integers}{\mathbb Z}

\newcommand{\abs}[1]{\vert#1\vert}

\newcommand{\setdiam}{\delta}
\newcommand{\genericbox}{{\mathcal A}}
\newcommand{\genericset}{{\mathcal E}}
\newcommand{\contourset}{{\mathcal F}}
\newcommand{\imcgfargv}{{\bolds\xi}}
\newcommand{\imcgfarg}{\xi}
\newcommand{\indexset}{{\mathbb J}}
\newcommand{\cgfdomain}{{\bolds\Theta}}
\newcommand{\bigset}{{\mathcal X}}
\newcommand{\varfcn}[1]{V_{#1}}
\newcommand{\latticedim}{d_0}
\newcommand{\fuzzy}{{\mathcal H}}
\newcommand{\transpose}{\top}
\newcommand{\shortdx}{d \bx}
\newcommand{\bp}{{\mathbf p}}
\newcommand{\by}{{\mathbf y}}

\makeatother

\begin{document}
\begin{frontmatter}

\title{Saddlepoint approximations for likelihood ratio like statistics
with applications to~permutation~tests}
\runtitle{Saddlepoint approximations for permutation tests}

\begin{aug}
\author[A]{\fnms{John} \snm{Kolassa}\corref{}\thanksref{t1}\ead[label=e1]{kolassa@stat.rutgers.edu}}
\and
\author[B]{\fnms{John} \snm{Robinson}\thanksref{t2}\ead[label=e2]{john.robinson@sydney.edu.au}}
\runauthor{J. Kolassa and J. Robinson}
\affiliation{Rutgers University and University of Sydney}
\address[A]{Department of Statistics\\
Rutgers University\\
110 Frelinghuysen Rd\\
Piscataway, New Jersey 08854-8019\\
USA\\
\printead{e1}}
\address[B]{School of Mathematics and Statistics\\
University of Sydney\\
NSW 2006\\
Australia\\
\printead{e2}} 
\end{aug}

\thankstext{t1}{Supported in part by NSF Grant DMS-09-06569.}
\thankstext{t2}{Supported in part by ARC DP0451722.}

\received{\smonth{4} \syear{2011}}
\revised{\smonth{8} \syear{2011}}

%
\begin{abstract}
We obtain two theorems extending the use of a saddlepoint approximation
to multiparameter problems for likelihood ratio-like statistics which
allow their use in permutation and rank tests and could be used in
bootstrap approximations. In the first, we show that in some cases when
no density exists, the integral of the formal saddlepoint density over
the set corresponding to large values of the likelihood ratio-like
statistic approximates the true probability with relative error of
order $1/n$. In the second, we give multivariate generalizations of the
Lugannani--Rice and Barndorff-Nielsen or $r^*$ formulas for the
approximations. These theorems are applied to obtain permutation tests
based on the likelihood ratio-like statistics for the $k$ sample and
the multivariate two-sample cases. Numerical examples are given to
illustrate the high degree of accuracy, and these statistics are
compared to the classical statistics in both cases.
\end{abstract}

%

\end{frontmatter}

\section{Introduction}
In parametric problems where distributions are specified exactly, the
likelihood ratio is
generally used for hypothesis testing whenever possible. In
multiparameter problems,
the distribution of twice the log likelihood ratio is approximated by a
chi-squared
distribution. Refinements of this approximation were obtained by
Barndorff-Nielsen \cite{BN86}
for parametric problems. In a nonparametric setting the empirical
exponential likelihood is described in Chapter 10 of \cite{DH} and
discussed in a number of references cited there. Saddlepoint
approximations for empirical exponential likelihood statistics based on
multiparameter $M$-estimates are given, for example, in \cite{RRY} and
for tests of means in \cite{R}, under the strong assumption that the
density of the $M$-estimate exists and has a saddlepoint approximation.
They used methods based on those of
\cite{BNC84} to obtain an approximation analogous to the Lugananni--Rice
approximation for the one-dimensional case.

It is the purpose of this paper to show that, under conditions
which will allow the application of the approximations in bootstrap,
permutation and rank statistics used for multiparameter cases, the
integral of the formal saddlepoint density approximation can be used
to give an approximation with relative error of order $n^{-1}$ to the
tail probability of a likelihood ratio-like statistic. This then permits
the approximation to be put in the Lugananni--Rice form as in \cite{RRY}
and also in a form analogous to the $r^*$ or Barndorff-Nielsen form given
in \cite{BN86} and \cite{FJR94} for the one-dimensional case. These
results are then applied to
two multiparameter nonparametric cases. We require the existence of a
moment generating function. This may be too strong an assumption in the
case of tests concerning means considered here, but robust versions of
these, as in \cite{RRY}, can be used to make the results widely applicable.

The statistic used is obtained by using the conjugate distribution
approach of \cite{C38} and is the log likelihood ratio in the
parametric case of exponential families. It can be written as a convex
function of $\bar{\mathbf X}$, the mean of $n$ independent random
variables. This statistic can be approximated to first order by a
quadratic form in the means~$\bar{\mathbf X}$. However, it does not
seem to be possible to approximate tail probabilities for quadratic
forms with relative errors of order $n^{-1}$, as are obtained for our
statistic. Cram\'{e}r large deviation results for the case of quadratic
forms in multivariate means were obtained by \cite{Osipov} and a number
of earlier authors cited in that paper, but the relative errors for the
approximation to the probability of the statistic, a random variable of
order $1/n$, exceeding $\lambda$ is of order $\sqrt{n\lambda}n^{-1/4}$.
So the relative error is at best of order $n^{-1/4}$. The same problem
arises in the case of an empirical likelihood statistic, where we know
of no saddlepoint approximation.

In the next section we introduce the notation and assumptions necessary
to obtain the likelihood ratio-like statistic, tail probabilities of
which can be used for hypothesis testing in multivariate nonparametric
settings. We reduce certain conditional cases given lattice variables
to a more convenient notation and state the main result in a theorem
showing that tail probabilities for the statistic can be approximated,
to relative order $n^{-1}$, by an integral of a formal saddlepoint
density. We then state and prove a theorem giving the integrals in
forms like those of Lugananni--Rice and Barndorff-Nielsen in the
one-dimensional case. In Section~\ref{sec3} we consider two examples of
permutation tests, for the $k$-sample problem and for a two sample
multivariate permutation test, using the results of the previous
section to obtain explicit formulas for test statistics and for the
approximations of the tail probabilities of these statistics under
permutations. We then present numerical examples illustrating the
accuracy of the approximations and comparing results to those obtained
using the standard sum of squares test statistics for the $k$-sample
permutation and rank tests and the Mahalanobis $D^2$ test for the
2-sample multivariate test. In the final section we give the proof of
the main result.

\section{Notation and main result}\label{sec2}

For a sample of size $n$ with mean vector $\bar\bx$ from a parametric
canonical exponential family with density $f_\cgfargv(x)=\exp(\cgfargv
^\top\bx-\kappa(\cgfargv))g(\bx)$, the maximum likelihood estimate of
$\cgfargv$ is $\hat\cgfargv$, the solution of $\kappa'(\cgfargv)=\bar\bx
$, and, taking $\kappa'(\mathbf{0})=\mathbf{0}$, the log likelihood
ratio statistic is $\Lambda(\bar\bx)=\hat\cgfargv{}^\top\bar\bx-\kappa
(\hat\cgfargv)$. This is used to test the hypothesis that $\cgfargv
=\mathbf{0}$, or equivalently, that $\kappa'(\cgfargv)=\mathbf{0}$. For
the nonparametric case an empirical exponential family is taken, and
it is shown, for example, in \cite{RRY}, page 1163, that the empirical
exponential likelihood ratio statistic for a test that the expectation
is zero is $\Lambda(\bar\bx)=-{\mathbf{\beta}}_0^\top\bar\bx+\kappa
_n(\mathbf{\beta}_0)$, where $\kappa_n(\mathbf{\beta})=\log[\sum
_{i=1}^n\exp(\mathbf{\beta}^\top\bx_i)]/n$ and $\mathbf{\beta}_0$ is
the solution of $\kappa'_n(\mathbf{\beta})=\mathbf{0}$. In \cite{RRY} a
bootstrap approximation can be based on the statistic $\Lambda(\bar\bx
^*)=\hat\cgfargv^\top\bar\bx^*-\kappa_n(\mathbf{\beta}_0+\hat\cgfargv)+
\kappa_n(\mathbf{\beta}_0)$, where the bootstrap is taken from the
tilted\vspace*{1pt} empirical distribution $\hat F_0(\bx)=\sum_{i=1}^n\exp(\mathbf
{\beta}_0^\top\bx_i-\kappa_n(\mathbf{\beta}_0))I\{\bx_i\leq\bx\}/n$. A
saddlepoint approximation to this bootstrap is given, but it is noted
that the relative errors of this approximation could not be proven from
the theorem of that paper. The theorems of this section permit this
proof. We use an analogous approach to give the likelihood ratio-like
statistics for the two permutation test examples in the next section.

Consider independent $d$-dimensional
random vectors $\bX_1,\ldots,\bX_n$, with the first
$d_0$ components $X_{1j},\ldots,X_{\latticedim j}$ confined to
a lattice with unit spacings, for $\latticedim<d$, and with the
average cumulant generating function
%
%
\begin{equation}\label{cgfdef}
\kappa(\cgfargv)=n^{-1}\log(Ee^{\cgfargv^\top\bS_n})=n^{-1}\sum
_{i=1}^n\log(Ee^{\cgfargv^\top\bX_i}),
\end{equation}
where $\bS_n=\bX_1+\cdots+\bX_n$. For some $\bx$ we can define
%
%
\begin{equation}\label{lambdadef}
\Lambda(\bx)=\hat{\cgfargv}{}^\top\bx-\kappa(\hat\cgfargv)
\end{equation}
for $\hat\cgfargv$ satisfying
%
%
\begin{equation}\label{spdef}
\kappa'(\hat{\cgfargv})=\bx
\end{equation}
and
%
%
\begin{equation}\label{rdef}
r(\bx)=e^{-n\Lambda(\bx)}
(2\pi n)^{-d_0/2}
(2\pi/n)^{-d_1/2}
|\varfcn{\hat{\cgfargv}}|^{-1/2}.
\end{equation}
In the case when the last $d_1$ components of $\bX_1,\ldots,\bX_n$ have
densities, this is the saddlepoint density approximation for $\bar\bX
=\bS_n/n$, obtained in the case of identically distributed random
vectors in \cite{BR65}. In many cases when these last components lack a
density, the theorem below will imply that their distribution may be
well approximated by a continuous distribution.

Let $\mu$ denote the distribution of $\bar\bX=\bS_n/n$, let
$\cgfdomain^*=\{\cgfargv\dvtx\kappa(\cgfargv)<\infty\}$, and let
$\mu_\cgfargv(d\byr)=\exp(-n(\kappa(\cgfargv)-\cgfargv^\transpose
\byr)) \mu(d\byr)$ define the distribution of $\bar\bX_\cgfargv$, the
mean of $\bX_{1\cgfargv},\ldots,\bX_{n\cgfargv}$, the associated
independent random vectors. These conjugate distributions,
first\vspace*{1pt} introduced in \cite{C38}, permit us to consider
large deviations. Let $V_\cgfargv=\kappa''(\cgfargv)$, and, taking
$\Vert\bx\Vert=(\bx^\top\bx)^{1/2}$, let
\[
\eta_j(\cgfargv)=n^{-1}\sum_{i=1}^n E\bigl[\Vert V_\cgfargv^{-1/2}
(\bX_{i\cgfargv}-E[\bX_{i\cgfargv}])\Vert^j\bigr].
\]
Let
%
%
\begin{eqnarray}\label{qdef}
q_{\cgfargv}(T)&=&\sup\bigl\{\bigl|e^{\kappa(\cgfargv+i\imcgfarg)-
\kappa(\cgfargv)}\bigr|
\dvtx\Vert\varfcn{\cgfargv}^{1/2}\imcgfargv\Vert>(3/4)\eta_3(\cgfargv
)^{-1},\nonumber\\[-8pt]\\[-8pt]
&& \hspace*{51.5pt}\abs{\imcgfarg_i}<\pi\mbox{ for }i\leq
d_0,\abs{\imcgfarg_i}<T,i>d_0\bigr\}.\nonumber
\end{eqnarray}

We consider the following conditions, essentially from \cite{RHHQ},
where, throughout, $c$ and $C$ are generic positive constants, and
$|A|$ denotes the determinant of a square matrix $A$. The complexity of
these conditions is due to the fact that we need to consider
conditional distributions of independent, but not identically
distributed, random variables.
\begin{itemize}
\item{(A1)} There is a compact subset, $\cgfdomain$, of the interior of
$\cgfdomain^*$,
with $\bmz$ in the interior of $\cgfdomain$.
\item{(A2)} $|V_\cgfargv|>c>0$ for $\cgfargv\in\cgfdomain$.
\item{(A3)} $\eta_j(\cgfargv)<C$ for $j=1,\ldots,5$ and $\cgfargv\in
\cgfdomain$.\vspace*{2pt}
\item{(A4)} $n^{2d_1+2} q_\cgfargv(n^{-2})<C$.
\end{itemize}
Here the first condition asserts that there is an open neighborhood of
the origin where the cumulative generating function exists. The second
condition bounds the average variance of the associated random
variables away from zero, and the third gives upper bounds the first 5
standardized moments in this neighborhood. The fourth condition is a
smoothness condition introduced first for the univariate case in \cite
{ABVZ} and which is sufficient to allow Edgeworth expansions for many
statistics based on ranks and applications to bootstrap and permutation
statistics when the original observations are from a continuous distribution.

Let $\bigset=\kappa'(\cgfdomain)$; then we are able to obtain
equations (\ref{lambdadef}), (\ref{spdef}) and (\ref{rdef}) for $\bx\in
\bigset$. Also, if $d_0>0$, let
$\Lambda_0(\bx_0)=\hat\cgfargv_0^\transpose\bx_0-\kappa_0(\hat\cgfargv_0)$
for $\hat\cgfargv_0$ satisfying \mbox{$\kappa'_0(\hat\cgfargv_0)=\bx_0$}, where
the subscript $0$ denotes a reduction to the first $d_0$ elements of
the $d$-vectors, and
we will use the subscript $1$ to denote the last $d_1$ elements.
If $r_0(\bx_0)=(2\pi n)^{-d_0/2}
|V_{\hat\cgfargv_0}|^{-1/2}\exp(-n\Lambda_0(\bx_0))$, then from \cite
{BR65}, we have
$\Prob(\bar\bX_0=\bx_0)=\mu_0(\bx_0)=r_0(\bx_0)(1+O(1/n))$.
For $d_0>0$, we will consider
$\bx^\transpose=(\bx_0^\transpose,\bx_1^\transpose)$ and replace
$r(\bx)$ by
%
%
\begin{equation}
r(\bx_1|\bx_0)=r(\bx)/r_0(\bx_0)=
{\frac{|V_{\hat\cgfargv_0}|^{1/2}e^{-n(\Lambda(\bx)-\Lambda_0(\bx
_0))}}{(2\pi/n)^{d_1/2} |V_{\hat\cgfargv}|^{1/2}}},
\end{equation}
and replace $\mu$ by the distribution of $\bar\bX_1$ conditional on
$\bar\bX_0=\bar\bx_0$,
so that we consider conditional probabilities of $\bar\bX_1$ given
$\bar\bX_0=\bx_0$, associated with sets of form
\[
\contourset=\{\bx_1\dvtx\Lambda(\bx)-\Lambda_0(\bx_0)\geq\lambda, \bx
^\transpose
=(\bx_0^\transpose,\bx_1^\transpose)\}.
\]

The main result is the following theorem, whose proof is deferred to a~later section.
%
%
\begin{theorem}\label{th1}
Under conditions \textup{(A1)--(A4)},
%
%
\begin{equation}\label{mainres}
\biggl|\mu(\contourset)-\int_\contourset r(\bx_1|\bx_0)\,d\bx_1\biggr|
=\biggl[\int_\contourset r(\bx_1|\bx_0)\,d\bx_1\biggr] O(1/n).
\end{equation}
\end{theorem}

Note that if the nonlattice subvectors of $\bX_1,\ldots,\bX_n$ have
densities, the variables are identically distributed and (A1) and (A2) hold,
then the theorem follows from Theorem 1 of \cite{BR65}.

The following theorem is a corollary whose derivation we include here.
This is the form that will be used in examples.
%
%
\begin{theorem}\label{rry}
Under the conditions of Theorem \ref{th1}, if $u=\sqrt{2\lambda}$,
%
%
\begin{equation}\label{hi}
\int_\contourset r(\bx_1|\bx_0)\,\shortdx_1=\bar Q_{d_1
}(nu^2)[1+O(1/n)]
+\frac{ c_n}{n}u^{d_1} e^{-nu^2/2} \frac{G(u)-1}{u^2}
\end{equation}
and
%
%
\begin{equation}\label{bn}
\int_{\contourset} r(\bx_1|\bx_0)\,\shortdx_1=
\bar Q_{d_1 }(nu^{*2})[1+O(1/n)],
\end{equation}
where $\bar Q_d(x)=\Prob(\chi_d^2\geq x)$,
%
%
\begin{eqnarray}\label{ust}
u^*&=&u-\log(G(u))/nu,\\
\label{cn}
c_n&=&\frac{n^{d_1 /2}}{2^{d_1 /2-1}\Gamma(d_1
/2)},\\
\label{deltadef}
\delta(\sqrt{2\lambda},s)&=&{
\frac{\Gamma(d_1 /2)|V_{\hat\cgfargv_0}|^{1/2}|V_{\hat\cgfargv
}|^{-1/2}|V_0|^{1/2} r^{d_1 -1}}
{2\pi^{d_1 /2}u^{d_1 -2}|\bs^\transpose V_{0}^{1/2}\hat
\cgfargv_1|}
},\\
\label{gg}
G(u)&=&\int_{S_{d_1 }} \delta(u,s)\,ds
\end{eqnarray}
for $S_{d_1 }$ the $d_1 $-dimensional unit sphere
centered at
zero, and where, for each $\bs\in S_{d_1 }$, $r$ is chosen so
$\Lambda(\bx_0,r \bs)-\Lambda_0(\bx_0)=\lambda$ and $V_0^{-1}=[\kappa
''(\mathbf{0})^{-1}]_{11}$, with the final subscripts denoting the lower
right $d_1\times d_1$ submatrix.
\end{theorem}
\begin{pf}
The derivation of (\ref{hi}), given Theorem \ref{th1}, is given in
\cite{RRY}. To get~(\ref{bn}), we use a related method. After making
the transformations $\byr=V_{\mathbf{0}}^{-1/2}\bx_1 $, $\byr\to(r,\bs)$
and $(r,\bs)\to(u,\bs)$, where the first is the polar transformation
with $\|\bx\|=r$ and $\bs\in S_d$, the unit sphere in $d$-dimensions,\vadjust{\goodbreak}
and the second has $u=\sqrt{2 (\Lambda (\bx_0,r \bs)-\Lambda_0(\bx_0))
}$, we have
\begin{eqnarray*}
\int_\contourset r(\bx_1|\bx_0) \,d\bx_1
&=&c_n \int_{u}^\infty v^{d-1}e^{-n v^2/2}G(v)\,d v\\[-3pt]
&=&c_n \int_{u}^\infty v^{d-1}e^{-n(v-\log G(v)/n
v)^2/2}\,dv\bigl(1+O(1/n)\bigr).\vspace*{-2pt}
\end{eqnarray*}
Then make the transformation $v^*=v-\log G(v)/n v$.
The final equality follows since $G(v)=1+v^2k(v)$ and $G'(v)=vk^*(v)$,
where $k(v)$ and $k^*(v)$ are bounded as shown in \cite{RRY}.\vspace*{-2pt}
\end{pf}
\begin{Remark*}
The integral (\ref{gg}) can
be approximated by a Monte Carlo method, for example, by approximating
$\int_{S_d}h(s)\,ds$ as
\[
\frac{2\pi^{d/2}}{\Gamma{(d/2)}}\frac{1}{M}\sum_{\ell=1}^M h(U_\ell),\vspace*{-2pt}
\]
where $U_1,\ldots,U_M$ are i.i.d. uniformly distributed on $S_d$. Here
the number of replicates in the Monte Carlo simulation can be small
with little loss of accuracy. We discuss this in the examples where it
was found that $M=10$ was sufficient. It would be possible to use a
method such as that in \cite{Genz} to get a numerical approximation to
the integral, but the Monte Carlo method is much simpler to use and
easily gives the required accuracy.\vspace*{-3pt}
\end{Remark*}

\section{Two examples of permutation tests}\label{sec3}
We consider a $k$ sample permutation test in a one-way design and a
multivariate two-sample permutation test.
In both cases we consider hypotheses that the populations of
random variables or vectors are exchangeable. In the first case the
observations are
generated either by sampling $n_1,\ldots,n_k$ independent random variables
from distributions $F_1,\ldots,F_k$, and we test $H_0\dvtx F_1=\cdots=F_k$,
or they are generated from an experiment in which $k$ treatments are
allocated at random to groups of sizes $n_1,\ldots,n_k$, and we test
$H_0$: treatments have equal effects. We choose a statistic suitable
for testing with respect to differences in means. The standard choices
of test statistic are the $F$-statistic from the analysis of variance or,
for a
nonparametric test based on ranks, the Kruskal--Wallis statistic. In
the second case the observations are generated by sampling from two
populations of $l$-dimensional random vectors, and we test for equality
of the distributions, or they are generated by experimental
randomization, and we test for equality of two treatments. Here the
test statistic arising from an assumption of multivariate normality is
the Mahalanobis $D^2$ test.\vspace*{-3pt}

\subsection{Permutation tests for $k$ samples}
Suppose that $a_1,\ldots,a_N$ are the elements of a finite population,
such that $\sum_{m=1}^N a_m=0$ and $\sum_{m=1}^N a_j^2=N$.
Let~$n_1,\allowbreak\ldots,n_k$ be integers, such that $N=\sum_{i=1}^k n_i$.
Suppose that $R_1,\ldots,R_N$ is an equiprobable random permutation of
$1,\ldots,N$. Let $X_{ij}=a_{R_{n_1+\cdots+n_{i-1}+j}}$, and let $\bar
X_i=\sum_{j=1}^{n_i} X_{i j}/n_i$.\vadjust{\goodbreak}

For $i=1,\ldots,k-1$, let $\be_i$ have $k-1$ components,
with component $i$ equal to $1$, and other components zero. Let
$\bI_m,  m=1,\ldots,N$ be independent and identically distributed
random vectors with $\Prob(\bI_m=\be_i)=n_i/N=p_i$
for $i<k$ and $\Prob(\bI_m=\bmz)=n_k/N=p_k$. Let
$\bS^\top=(\sum_{m=1}^N\bI_m^\top,\sum_{m=1}^Na_m\bI_m^\top)=(\bS_0^\top
,\bS_1^\top)$.
We have
\[
\Prob(n\bar\bX\leq\bx)=\Prob(\bS_1\leq
\bx|\bS_0=N\bp),
\]
where $\bar\bX,  \bx,  \bp$ are $k-1$ vectors
corresponding to the first $k-1$ samples.
Under~$H_0$, the cumulant generating
function of $\bS$ is
\[
N \kappa(\bolds\tau_0,\bolds\tau_1)=\log E
e^{\sum_{m=1}^N(\bolds\tau_0^\top\bI_m+\bolds\tau_1^\top\bI_ma_m)}
=\sum_{m=1}^N\log\Biggl(p_k+\sum_{i=1}^{k-1}p_i
e^{\tau_{0i}+\tau_{1i}a_m}\Biggr).
\]
Let $(\hat{\bolds{\tau}}_0^\top,\hat{\bolds{\tau}}_1^\top)$
be the solution of
\[
\kappa'(\bolds\tau_0,\bolds{\tau}_1)=(\bp,\bx).
\]

Let $B=\{\bx\dvtx\Lambda(\bx)\geq u^2/2\}$, where $\Lambda(\bx)=\hat{\bolds
\tau
}_0^\top
\bp+\hat{\bolds\tau}_1^\top\bx-\kappa(\hat{\bolds\tau}_0,\hat{\bolds\tau}_1)$,
and note that $\kappa'(0,\bmz)=(p,\bmz)$ and $\kappa(0,\bmz)=0$.
Now from Theorem \ref{th1}, if $q_{\bolds\tau}(n^{-2})=O(n^{-2k})$,
\[
\Prob\bigl(\Lambda(\bar\bX)\geq u^2/2\bigr)=\int_B
r(\bx|\bp)\,d\bx\bigl(1+O(1/N)\bigr),
\]
where
\[
r(\bx|\bp)
=(2\pi/N)^{-(k-1)/2}|\kappa_{00}(0,\bmz)|^{1/2}|\kappa''(\hat{\bolds\tau
}_0,\hat{\bolds\tau}_1)|^{-1/2}e^{-N\Lambda(\bx)}.
\]
Then from Theorem \ref{rry},
$G(u)$ is given in (\ref{deltadef}) and (\ref{gg}) with $d_0=d_1=k-1$.
Now we can use (\ref{hi}) and (\ref{bn}) to get the two approximations.

\subsection{Numerical results for $k$-sample test} Consider first the
rank test based
on the statistic $\Lambda(\bar\bX)$ where $ a_1=1,
\ldots, a_N=N$, with $N=20$ for 4 groups of size~5; in
the standard case the Kruskal--Wallis test would be used. The following
table gives the results of tail probabilities from a Monte Carlo
simulation of $\Lambda(\bar\bX)$ (MC~$\Lambda$) and of the
Kruskal--Wallis statistic (MC K--W)
using 100,000 permutations, the chi-squared approximation ($\chi^2_3$)
and the
saddlepoint approximations using (\ref{hi}) (SP LR $\Lambda$) and (\ref
{bn}) (SP BN $\Lambda$), using $M=1000$
Monte Carlo samples from $S_3$. Inspection of the table comparing
the saddlepoint Lugananni--Rice and Barndorff-Nielsen approximations
with the Monte Carlo
approximation for $\Lambda$
shows the considerable accuracy of these approximations throughout the range.
The chi square approximations to the distribution of the
Kruskal--Wallis statistic does not have this degree
of accuracy. We note that good approximations for the saddlepoint
approximations are achieved by $M$ as small as 10.
We obtained the standard deviation of individual random values of the
integrand and noted that for Table \ref{table1} this was 0.003 for $\hat u=0.6$ and
0.0007 for $\hat u=0.9$, indicating that $M=10$ gives sufficient accuracy
in this example.\vadjust{\goodbreak}

%
%
\begin{table}
\caption{The $4$-sample rank tests with $n_i=5$}
\label{table1}
\begin{tabular*}{\tablewidth}{@{\extracolsep{\fill}}lccccccc@{}}
\hline
$\bolds{\hat u}$
& \textbf{0.3} & \textbf{0.4} & \textbf{0.5} & \textbf{0.6}
& \textbf{0.7} & \textbf{0.8} & \textbf{0.9}\\
\hline
MC $\Lambda$& 0.6758 &0.4328 &0.2365& 0.1087& 0.0423&0.0142&0.0041\\
MC K--W & 0.6583& 0.4027& 0.1921&0.0652& 0.0135&0.0012&0.0000\\
$\chi^2_3$& 0.6149& 0.3618& 0.1718 &0.0658& 0.0203& 0.0051 & 0.0010\\
SP LR $\Lambda$& 0.6811& 0.4446 &0.2454& 0.1151& 0.0464&0.0164&0.0052\\
SP BN $\Lambda$& 0.6753& 0.4380& 0.2387& 0.1101& 0.0434&0.0148&0.0045\\
\hline
\end{tabular*}
\end{table}

Also consider the permutation test based on a single sample of 40 in 4
groups of 10 from an exponential distribution, comparing as above each
of the saddlepoint approximations with the Monte Carlo approximations
in this case and with the standard test based on the sum of squares from
an analysis of variance. The same pattern of accuracy as reported above
is apparent from inspection of Table \ref{table2}.

%
%
\begin{table}[b]
\caption{The $4$-sample permutation tests with exponentially
distributed errors and $n_i=10$}
\label{table2}
\begin{tabular*}{\tablewidth}{@{\extracolsep{\fill}}lccccccc@{}}
\hline
$\bolds{\hat u}$ & \textbf{0.2} & \textbf{0.3} &
\textbf{0.4} & \textbf{0.5} & \textbf{0.6} & \textbf{0.7} & \textbf{0.8}\\
\hline
MC $\Lambda$ & 0.9456 &0.6837& 0.3434& 0.1160& 0.0275& 0.0043& 0.0004\\
MC ANOV & 0.9455 &0.6784& 0.3273& 0.0971& 0.0164& 0.0015&0.0004\\
$\chi^2_3$& 0.9402& 0.6594&0.3080& 0.0937& 0.0186& 0.0024 &0.0002 \\[1pt]
SP LR $\Lambda$& 0.9491 &0.6888& 0.3456& 0.1174& 0.0272& 0.0043
&0.0004\\
SP BN $\Lambda$& 0.9486& 0.6877 &0.3441 &0.1164& 0.0268 &0.0042&0.0004\\
\hline
\end{tabular*}
\end{table}

\subsection{A two-sample multivariate permutation test}
Let $\mathbf{a}_1,\ldots,\mathbf{a}_N$ be $l$-vectors regarded as
elements of a finite population such that $\sum_{i=1}^N\mathbf{a}_i=0$
and $\sum_{i=1}^N\mathbf{a}_i\mathbf{a}_i^T=NI$. Let $R_1,\ldots,R_N$
be obtained by an equiprobable random permutation of $1,\ldots,N$, let
$\mathbf{X}_j=\mathbf{a}_{R_j}$, $j=1,\ldots,N$ and\vspace*{1pt} let $\bar{\mathbf
X}_1=\sum_{j=1}^n\mathbf{X}_j/n$ for $n=Np$ with $0<p<1$.
Let $I_1,\ldots,I_N$ be i.i.d. Bernoulli variables with $EI_1=p$. If
$\mathbf{S}^T=(S_0,\mathbf{S}_1^T)$ with $S_0=\sum_{i=1}^NI_i$ and
$\mathbf{S}_1=\sum_{i=1}^N\mathbf{a}_iI_i$, then for any Borel set
$\mathcal F$,
%
%
\begin{equation}\label{mvtsp}
P(\bar{\mathbf X}\in\mathcal F)=P(\mathbf{S}_1/N\in\mathcal
F|S_0/N=p).
\end{equation}
Let $\bolds\tau^\top=(\tau_0,\bolds\tau_1^\top)$ with $\tau_0\in\Re$
and $\bolds
\tau_1\in\Re^d$ and let
\begin{eqnarray*}
\kappa(\bolds\tau)&=&N^{-1}\log E \exp(\tau_0S_0+\bolds\tau_1^\top
\mathbf
{S}_1)\\
&=&N^{-1}\sum_{i=1}^N\log(q+pe^{\tau_0+\bolds\tau_1^\top\mathbf{a}_i}).
\end{eqnarray*}
Let $\hat\cgfargv$ be the solution of $\kappa'(\bolds\tau)=(p,\mathbf
{x}^\top)^\top$, and let $\Lambda(p,\mathbf{x})=\hat\tau_0p+\hat\cgfargv
{}^\top_1\mathbf{x}-\kappa(\hat\cgfargv)$. Consider sets $\mathcal F=\{
\mathbf{x}\dvtx\Lambda(p,\mathbf{x})\geq\lambda\}$. Then from Theorem \ref
{rry}, we can approximate (\ref{mvtsp}) by (\ref{hi}) or (\ref{bn})
where $G(u)$ is given by (\ref{deltadef}) and (\ref{gg}) with $d_0=1$
and $d_1=l$.

\subsection{Numerical results for two-sample test}
Consider the test based on two samples of size 40 from a 3-variate
exponential distribution with mean~1 and covariance matrix $I$. After
standardizing the combined sample we consider tests based on the
statistic $\Lambda(\bar{\mathbf X})$ or $\bar{\mathbf X} \bar{\mathbf
X}^T$, equivalent to the usual normal theory based statistic. We
calculate the tail probabilities based on Theorem \ref{rry} in this case and
Monte Carlo approximations to the permutation tests based on 10,000
random permutations. Table \ref{table3} demonstrates the accuracy of the two
saddlepoint approximations throughout the range. It also shows that the
chi-squared approximation is not satisfactory either for $\Lambda$ or
for the classical quadratic form statistic. However, while we have
accurate tail probability approximations for the new statistic, such
approximations are not available for the classical quadratic form.

%
%
\begin{table}
\caption{The $3$-dimensional two sample parmutation test}
\label{table3}
\begin{tabular*}{\tablewidth}{@{\extracolsep{\fill}}lccccc@{}}
\hline
$\bolds{\hat u}$& \textbf{0.3} & \textbf{0.4} & \textbf{0.5}
& \textbf{0.6} & \textbf{0.7}\\
\hline
MC $\Lambda$ & 0.3543& 0.1249& 0.0276& 0.0041& 0.0006\\
$\chi_3^2$ & 0.3080& 0.0937& 0.0186 & 0.0024& 0.0002\\
SP LR $\Lambda$ & 0.3528 & 0.1207 & 0.0282& 0.0045 & 0.0005\\
SP BN $\Lambda$ & 0.3507& 0.1194& 0.0278 & 0.0043& 0.0005\\
Quadratic & 0.3325 & 0.0939& 0.0135 & 0.0004& 0.0001\\
\hline
\end{tabular*}

\end{table}

\section{Proofs of the main results}

For notational convenience we will restrict attention to the case
$d_0=0$, as details of the case conditional on lattice variables follow
in a straightforward manner.
The following theorem is a~simplified version of Theorem 1 of
\cite{RHHQ}, taking $s=5$, $d_0=0$ and $\mathcal A$ as a~$d$-dimensional cube in $\bigset$ with center $\ba$ and side
$\delta=n^{-1}$.
As in (1.10) of~\cite{RHHQ}, let
\[
e_2(\byr,\bmg_\tau)=\bigl(1+Q_1(\byr^*)+Q_2(\byr^*)\bigr)
(2\pi/n)^{-d/2}|V_\cgfargv|^{-1/2}e^{-\byr^{*\top}\by^*/2}
\]
with $\byr^*=n^{1/2} V_{\cgfargv}^{-1/2}(\byr-\kappa'(\cgfargv))$, be
the formal Edgeworth expansion of order~2 for $\bar\bX_{\cgfargv}=\sum_{i=1}^n
\bX_{i\cgfargv}/n$, and let
\[
e_2\bigl(\cgfargv,{\mathcal E},\bx-\kappa(\cgfargv)\bigr)=\int
_{\mathcal E}e^{n\cgfargv^\transpose(\bx-\byr)}e_2(\byr,\bmg_\tau)\,d\byr.
\]
The terms $Q_1$ and $Q_2$ are given explicitly in (1.11) of \cite
{RHHQ}, and
are terms of order $n^{-1/2}$ and $n^{-1}$, respectively.
%
%
\begin{theorem}\label{rhhq}
For any set $\genericset\subset{\mathcal A}$ and $\varepsilon>0$, take
$\genericset_{\varepsilon}=
\{\bz\dvtx\exists\byr\in\genericset,\Vert\bz-\byr\Vert<\varepsilon\}$.
Choose $\varepsilon\in(0,c/n^{2})$, and let $T=1/\varepsilon$.
For $\bx\in\genericset\subset\bigset$,
\begin{eqnarray*}
\bigl|\mu(\genericset)-e^{-n(\cgfargv^\top\bx-\kappa(\cgfargv))}
e_{2}\bigl(\cgfargv,\genericset,\bx-m(\cgfargv)\bigr)
\bigr|
\leq
e^{-n(\cgfargv^\top\bx-\kappa(\cgfargv))}
|\varfcn{\cgfargv}|^{-1/2}R
\end{eqnarray*}
for
\[
R=C\bigl[\Vol(\genericset_{2\varepsilon}) \bigl(\eta_{5}(\cgfargv)n^{-3/2}
+|\varfcn{\cgfargv}|^{1/2}n^{1/2}T^dq_{\cgfargv}(T) \bigr)+
\Vol(\genericset_{2\varepsilon}-\genericset_{-2\varepsilon}) \bigr].
\]
\end{theorem}

Note that this follows, since
\[
\hat\chi_{\cgfargv,\genericset_{2\varepsilon}}(\bmz)=
\int_{\genericset_{2\varepsilon}}e^{n\cgfargv^\transpose(\bu-\ba)}\,d u<C\Vol
(\genericset_{2\varepsilon})
\]
as $\genericset\subset{\mathcal A}$ implies that
$\Vert\bu-\ba\Vert< c(\delta+2\varepsilon)<c n^{-1}$.

We give a preliminary lemma before proceeding to the proof of Theorem~\ref{th1}, using the notation $\kappa(\bolds\tau(\bx))=\bx$ and $\Lambda(\bx
)=\bolds
\tau(\bx)^\top\bx-\kappa(\bolds\tau(\bx))$ for $\bx\in\bigset$.
%
%
\begin{lemma}\label{prelem}
For $\bx\in\genericset\subset{\mathcal A}\subset\bigset$,
%
%
\begin{equation}\label{diff1}
\int_\genericset r(\by)\,d\by-e^{n(\kappa(\bolds\tau(\bx))-\bolds\tau(\bx
)^\top
\bx)}e_2(\bolds\tau,\genericset,0)=\int_\genericset r(\by)\,d\by O(1/n).
\end{equation}
\end{lemma}
\begin{pf} Ignoring for the moment the terms involving $Q_1$ and $Q_2$,
the left-hand side in (\ref{diff1}) is
%
%
\begin{equation}\label{diff11}
\int_\genericset r(\by)\biggl[1-\frac{e^{n(\Lambda(\by)-\Lambda(\bx
)-\bolds
\tau(\bx)^\top(\by-\bx)-(\by-\bx)^\top V_{\bolds\tau(\bx)}^{-1}(\by-\bx
)/2)}}{| V_{\bolds\tau(\bx)}|^{1/2}/| V_{\bolds\tau(\by)}|^{1/2}}
\biggr]\,d\by.
\end{equation}
Noting that $\|\by-\bx\|=O(1/n)$, and using a Taylor series expansion
about~$\bx$, we see that the exponent in (\ref{diff11}) is $O(1/n^2)$,
and the denominator is $1+O(1/n)$. So in (\ref{diff1}) the first term
on the left is as given by the expression on the right. Noting that
$Q_1(\mathbf0)=0$, we see that the term involving $Q_1$ is of the same
form. The proof is completed by noting that the term $Q_2$ is also of
this form.
\end{pf}
\begin{pf*}{Proof of Theorem \ref{th1}}
The proof will proceed by dividing $\bigset$ into small rectangles, applying
Theorem \ref{rhhq} on each of these rectangles, and summing the
results
in a manner similar to that of \cite{Osipov}.
For $\bj\in\integers^d_1$, let
$\genericbox_\bj=\{\bx\in\reals^d\dvtx
x_l\in((j_{l-d_0}-\half)\setdiam,(j_{l-d_0}+\half
)\setdiam]\}$, and
let $\genericset^\bj=\genericbox_\bj\cap\contourset$.
By the intermediate value theorem, on each $\genericset^\bj$, there is
an $\bx_\bj$ such that
\[
\int_{\genericset^\bj} r(\bx)\,d\bx=r(\bx(\bj))
\Vol(\genericset^\bj).\vadjust{\goodbreak}
\]
Note that $\contourset=\bigcup_{\bj\in\indexset}\genericset^\bj$ and
$\genericset^\bj$ are disjoint. Define $\hat\cgfargv_\bj$ so that
$\bx_\bj=\kappa'(\hat\cgfargv_\bj)$.
Write $\indexset=\{j\dvtx\Vol(\genericset^\bj)>0\}$. Then
\begin{eqnarray*}
\mu(\contourset)-\int_{\contourset} r(\bx)\,d\bx
&=&\sum_{\bj\in\indexset}[\mu(\genericset^\bj)-r(\bx_\bj)
\Vol(\genericset^\bj)]\\
&=&E_1+E_2,
\end{eqnarray*}
where
\[
E_1=\sum_{\bj\in\indexset}[\mu(\genericset^\bj)
-r(\bx_\bj)(2\pi/n)^{d/2}|V_{\hat\cgfargv_\bj}|^{1/2}
e_2(\hat\cgfargv_\bj,\genericset^\bj,0)
]
\]
and
\[
E_2=-\sum_{\bj\in\indexset}\biggl[ \int_{\genericset^\bj} r(\by) \,d\by-
e^{n\kappa(\hat\cgfargv_\bj)-\hat\cgfargv_\bj^\top\bx_\bj}e
_2(\hat\cgfargv_\bj,\genericset^\bj,0) \biggr].
\]
Using Lemma \ref{prelem} on each $\genericset^\bj$ and summing, we have
\[
E_2=\int_\contourset r(\by)\,d\by O(1/n).
\]

Now consider $E_1$. Apply Theorem \ref{rhhq} to each $\genericset^\bj$,
and sum to get
%
%
\begin{equation}\label{E1}
\vert E_1\vert\leq\sum_{\bj\in\indexset}r(\bx_{\bj})(R_{1\bj}+R_{2 \bj}),
\end{equation}
where
\[
R_{1 j}=C \Vol(\genericset_{2\varepsilon}^\bj)
[\eta_5(\cgfargv_\bj)n^{-1}+
|V_{\hat\cgfargv_\bj}|^{1/2} n^{d/2} T q_{\hat\cgfargv_\bj
}(T)]
\]
and
\[
R_{2 \bj}=\Vol(\genericset_{2\varepsilon}^\bj-\genericset_{-2\varepsilon}^\bj).
\]
The summation of these terms is complicated by the fact that the sets
are not disjoint and not all are subsets of $\contourset$. So introduce
sets $\fuzzy_\bj=\genericbox_\bj\cap\contourset_{2\varepsilon}$. Consider
the set $\fuzzy^*_\bj$, the union of $\fuzzy_\bj$ and the $3^d-1$ sets
formed by reflections of $\fuzzy_\bj$ in each of the lower-dimensional
faces of $\genericbox_\bj$. Then $\genericset_\bj\subset\fuzzy^*_\bj$ so
\[
\Vol(\genericset_{\bj 2\varepsilon})/\Vol(\fuzzy_{\bj})\leq3^d.
\]
So
\[
\sum_{\bj\in\indexset}r(\bx_\bj)R_{1\bj}\leq\sum_{\bj\in\indexset}r(\bx
_\bj)\Vol(\fuzzy_\bj)O(1/n)= \int_{\contourset_{2\varepsilon}} r(\by)\,d\by O(1/n).
\]
Also
\[
\Vol(\genericset^\bj_{2\varepsilon}-\genericset^\bj_{-2\varepsilon})/\Vol
(\fuzzy_{\bj})\leq C\varepsilon/\delta=O(1/n).
\]
Using this to bound the second sum on the right-hand side of (\ref{E1}) and the
previous bound for the first term gives
\[
|E_1|=\int_{\contourset_{2\varepsilon}}r(\by)\,d\by O(1/n).
\]

Note that for any $\bx$ such that $\Lambda(\bx)=\lambda$ and any $\bz\in
\contourset_{2\varepsilon}$,
\[
\Lambda(\bz)\geq\Lambda(\bx)-|(\bz-\bx)^\top\Lambda'(\bx)|\geq\lambda
-C\varepsilon.
\]
So the theorem follows by noting that
\[
\hspace*{93pt}
\int_{\contourset_{2\varepsilon}}r(\by)\,
d\by=\int_{\contourset}r(\by)\,d\by\bigl(1+O(1/n)\bigr).
\hspace*{81pt}\qed
\]
\noqed\end{pf*}



\printaddresses

\end{document}